\documentclass[11pt]{amsart}

\newtheorem{theorem}{Theorem}

\newtheorem{corollary}[theorem]{Corollary}

\newtheorem{lemma}[theorem]{Lemma}

\newtheorem{proposition}[theorem]{Proposition}
\newtheorem{remark}[theorem]{Remark}

\usepackage{hyperref}

\usepackage[active]{srcltx} 

\usepackage{stackrel}

\def\J#1#2#3{ \left\{ #1,#2,#3 \right\} }
\def\RR{{\mathbb{R}}}
\def\NN{{\mathbb{N}}}
\def\11{\textbf{$1$}}
\def\CC{{\mathbb{C}}}

\usepackage{enumerate}
\usepackage{amsmath}
\usepackage{amssymb}

\begin{document}

\title[A Kaplansky Theorem for JB*-triples]{A Kaplansky Theorem for JB*-triples}

\author[Fern\'{a}ndez-Polo]{Francisco J. Fern\'{a}ndez-Polo}
\email{pacopolo@ugr.es}
\address{Departamento de An{\'a}lisis Matem{\'a}tico, Facultad de
Ciencias, Universidad de Granada, 18071 Granada, Spain.}

\author[Garc{\' e}s]{Jorge J. Garc{\' e}s}
\email{jgarces@correo.ugr.es}
\address{Departamento de An{\'a}lisis Matem{\'a}tico, Facultad de
Ciencias, Universidad de Granada, 18071 Granada, Spain.}

\author[Peralta]{Antonio M. Peralta}
\email{aperalta@ugr.es}
\address{Departamento de An{\'a}lisis Matem{\'a}tico, Facultad de
Ciencias, Universidad de Granada, 18071 Granada, Spain.}

\thanks{Authors partially supported by
D.G.I. project no. MTM2008-02186, and Junta de Andaluc\'{\i}a
grants FQM0199 and FQM3737.}
\thanks{Published at Proceeding of the American Mathematical Society \url{ https://doi.org/https://doi.org/10.1090/S0002-9939-2012-11157-8}. This manuscript version is made available under the CC-BY-NC-ND 4.0 license \url{https://creativecommons.org/licenses/by-nc-nd/4.0/.} }

\date{}

\begin{abstract} Let $T:E\rightarrow F$ be a
non-necessarily continuous triple homomorphism from a (complex)
JB$^*$-triple (respectively, a (real) J$^*$B-triple) to a normed
Jordan triple. The following statements hold:\begin{enumerate}
\item $T$ has closed range whenever $T$ is continuous. \item $T$
is bounded below if and only if $T$ is a triple monomorphism.
 \end{enumerate} \vspace{0.5 mm} This result generalises classical theorems of I.
Kaplansky \cite{Kaps} and S.B. Cleveland \cite{CLev} in the
setting of C$^*$-algebras and of A. Bensebah \cite{Ben} and J.
P{\'e}rez, L. Rico and A. Rodr{\'\i}guez Palacios \cite{PeRiRo} in the
setting of JB$^*$-algebras.
\end{abstract}

\maketitle
 \thispagestyle{empty}

\section{Introduction}

A celebrated result of I. Kaplansky (cf. \cite[Theorem
6.2]{Kaps}) establishes that any algebra norm on a commutative
C$^*$-algebra dominates the C$^*$-norm. Subsequently, S.B.
Cleveland (see \cite{CLev}) generalised this result to the
noncommutative case by showing that every (non necessarily
complete nor continuous) algebra norm on a $C^*$-algebra generates
a topology stronger than the topology of the C$^*$-norm. In other
words, every non necessarily continuous monomorphism from a
C$^*$-algebra to an associative normed algebra is bounded
below. Alternative proofs to Cleveland's result were given by H.G.
Dales \cite{Dales} and A. Rodr{\'\i}guez Palacios \cite{Ro90} (see also \cite[Theorem 6.1.16]{Pal}).\smallskip

The arguments presented by A. Rodr{\'\i}guez Palacios in \cite{Ro90}
were adapted by A. Bensebah \cite{Ben} and J. P{\'e}rez, L. Rico and A. Rodr{\'\i}guez Palacios
\cite{PeRiRo} to extend Kaplansky theorem to the more general setting of JB$^*$-algebras.
The results established in \cite{Ben} and \cite{PeRiRo} show that every non necessarily continuous Jordan
monomorphism from a JB$^*$-algebra to a normed Jordan algebra is
bounded below. This result was latter re-proved by S. Hejazian and
A. Niknam in \cite{HejNik}.\smallskip

Every C$^*$-algebra, $A$, admits a triple product defined by
\begin{equation}\label{equ triple prod c*} \J abc := \frac12 (a
b^* c +c b^*a).\end{equation} Let us suppose that $\|.\|_2$ is
another (non necessarily complete nor continuous) norm on $A$
which makes continuous the triple product of $A$. It is natural to
ask whether this norm generates a topology stronger than the
topology of the C$^*$-norm.\smallskip

Every C$^*$-algebra, $A$, equipped with its C$^*$-norm and the
triple product defined in $(\ref{equ triple prod c*})$ can be
regarded as an element in the wider category of (complex)
JB$^*$-triples (see \S 2 for the detailed definitions). The
question posed in the above paragraph also makes sense in the
(larger) categories of (complex) JB$^*$-triples and real
J$^*$B-triples. In this setting, the problem can be reformulated
in the following terms:\smallskip

\textbf{Problem $\mathbf{(P)}$} \emph{Let $E$ be a (complex)
JB$^*$-triple or a (real) J$^*$B-triple whose norm is denoted by
$\|.\|$, and let $\|.\|_2$ be another (non necessarily complete
nor $\|.\|$-continuous) norm on the vector space $E$ which makes
continuous the triple product of $E$. Does $\|.\|_2$ generate a
topology stronger than the topology generated by the JB$^*$-triple
norm $\|.\|$?}\smallskip

\emph{Equivalently,  is every (non-necessarily continuous) triple
monomorphism $T$ from $E$ to a normed Jordan triple bounded
below?}
 \smallskip

Under the additional hypothesis of $T$ being $\|. \|$-continuous
(resp., $\|.\|_2$  being $\|. \|$-continuous), Problem $(P)$ was
solved by K. Bouhya and A. Fern{\'a}ndez in the case of (complex)
JB$^*$-triples \cite[Corollary 14]{BouFer}.\smallskip

In this paper we solve Problem $(P)$ without any additional
assumptions on the triple monomorphism $T$ (resp., on $\|.\|_2$).
When particularized to C$^*$-algebras, our main result shows that
every non necessarily continuous triple monomorphism from a real
or complex C$^*$-algebra to a normed Jordan triple is bounded
below. \smallskip

Section \S 2 is devoted to present the basic facts and definitions
needed in the paper. We shall also survey the results on the
property of minimality of norm topology in the setting of Banach
algebras and Jordan-Banach triples. We shall adapt the arguments
given by K. Bouhya and A. Fern{\'a}ndez \cite{BouFer}, to obtain their
result in  the setting of real J$^*$B-triples.\smallskip

In Section \S 3 we present our main results (Theorem \ref{t
kaplansky without conrinuity} and Corollary \ref{c
Yood-Cleveland}). This section contains a deep study of the
separating spaces associated with a triple homomorphism between
normed Jordan triples. Among the tools developed here, we remark a
main boundedness theorem type for Jordan Banach triples (see
Theorem \ref{t boundeness}), which is the Jordan triple version of
a classical result in the setting of Banach algebras due to W.G.
Bade and P.C. Curtis \cite{BaCur}.

\section{Minimality of norm topology for JB$^*$-triples}

A normed algebra $A$ has \emph{minimality of algebraic norm
topology (MOANT)} if any other (non-necessarily complete) algebra
norm dominated by the given norm yields an equivalent topology. It
is part of the folklore that $C^*$-algebras have MOANT (compare
\cite[Lemma 5.3]{CLev}).\smallskip

In this section, we study the minimality of norm topology in the
setting of normed Jordan triples. We recall that a complex (resp.,
real) \emph{normed Jordan triple} is a complex (resp., real)
normed space $E$ equipped with a non-trivial, continuous triple
product $$ E \times E \times E \rightarrow E$$
$$(x,y,z) \mapsto \J xyz $$
which is bilinear and symmetric in the outer variables and
conjugate linear (resp., linear) in the middle one satisfying the
so-called \emph{``Jordan Identity''}:
$$L(a,b) L(x,y) -  L(x,y) L(a,b) = L(L(a,b)x,y) - L(x,L(b,a)y),$$
for all $a,b,x,y$ in $E$, where $L(x,y) z := \J xyz$. If $E$ is
complete with respect to the norm (i.e. if $E$ is a Banach space),
then it is called a complex (resp., real) \emph{Jordan-Banach
triple}. Every normed Jordan triple can be completed in the usual
way to become a Jordan-Banach triple. Unless otherwise is
specified, the term ``normed Jordan triple'' (resp.,
``Jordan-Banach triple'') will always mean a real or complex
normed Jordan triple (resp., ``Jordan-Banach triple'').\smallskip

For each Jordan-Banach triple $E$, the constant $N(E)$ or
$N(E,\|.\|)$ will denote the supremum of the set $\{\|\J xyz \|:
\|x\|,\|y\|,\|z\|\leq 1\}$.\smallskip

A real (resp., complex) \emph{Jordan algebra} is a
(non-necessarily associative) algebra over the real (resp.,
complex) field whose product is abelian and satisfies $(a \circ
b)\circ a^2 = a\circ (b \circ a^2)$. A normed Jordan algebra is a
Jordan algebra $A$ equipped with a norm, $\|.\|$, satisfying $\|
a\circ b\| \leq \|a\| \ \|b\|$, $a,b$ in $A$. A \emph{Jordan
Banach algebra} is a normed Jordan algebra whose norm is
complete.\smallskip

Every real or complex associative Banach algebra (resp., Jordan
Banach algebra) is a real Jordan-Banach triple with respect to the
product $\J abc = \frac12 (a bc +cba)$ (resp., $\J abc = (a\circ
b) \circ c + (c\circ b) \circ a - (a\circ c) \circ b$).\smallskip

A JB$^*$-algebra is a complex Jordan Banach algebra $A$ equipped
with an algebra involution $^*$ satisfying that $\|\J a{a^*}a \| =
\| 2 (a\circ a^*) \circ a - a^2 \circ a^*\|= \|a\|^3$, $a$ in
$A$.\smallskip

Every JB$^*$-algebra has MOANT (compare \cite[Theorem
10]{PeRiRo}).\smallskip

We shall say that a normed Jordan triple $E$ has \emph{minimality
of triple norm topology (MOTNT)} if any other (non-necessarily
complete) triple norm dominated by the norm of $E$ defines an
equivalent topology.\smallskip

\begin{remark}
\label{r MOANT equivalent to MOTNT}{\rm  Let $A$ be a real or
complex associative normed algebra whose norm is denoted by
$\|.\|$. The symbol $A^{+}$ will stand for the normed Jordan
algebra $A$ equipped with the Jordan product $a\circ b = \frac12
(a b +ba)$ and the original norm. Let $\|.\|_1$ be a norm on the
space $A$. Since the Jordan product is $\|.\|_1$-continuous
whenever the associative product is, we deduce:
$$(A^{+},\|.\|) \hbox{ has MOANT } \Longrightarrow
(A,\|.\|) \hbox{ has MOANT}.$$ However, we do not know if the
reciprocal statement is, in general, true. By \cite[Proposition
3]{CaMorRo}, there exists an associative normed algebra
$\mathcal{B}$ such that there exists a norm $\|.\|_1$ on
$\mathcal{B}$ for which the Jordan product is continuous but the
associative product is discontinuous. In particular,
$(\mathcal{B}^{+},\|.\|_1)$ doesn't have MOANT.\smallskip

When $A$ is simple and has a unit, every norm on $A$ making the
Jordan product continuous also makes continuous the associative
product (compare \cite[Theorem 3]{CaMorRo}). Under this additional
hypothesis, we have $$(A^{+},\|.\|) \hbox{ has MOANT }
\Longleftrightarrow (A,\|.\|) \hbox{ has MOANT}.$$

Suppose that $J$ is a real or complex normed Jordan algebra, whose
norm is denoted by $\|.\|$. When $J$ is regarded as a real or
complex normed Jordan triple with respect to the product $\J abc =
(a\circ b) \circ c + (c\circ b) \circ a - (a\circ c) \circ b$,
every Jordan algebra norm on $J$ makes continuous the triple
product. Therefore $J$ has MOANT whenever it has MOTNT.\smallskip

When $J$ has a unit, the Jordan and the triple product of $J$ are
mutually determined, and hence $$(J,\|.\|) \hbox{ has MOANT }
\Longleftrightarrow (J,\|.\|) \hbox{ has MOTNT}.$$}
\end{remark}\smallskip

A \emph{(complex) JB$^*$-triple} is a complex Jordan Banach triple
${E}$ satisfying the following axioms: \begin{enumerate}[($JB^*
1$)] \item For each $a$ in ${E}$ the map $L(a,a)$ is an hermitian
operator on $E$ with non negative spectrum. \item  $\left\|
\{a,a,a\}\right\| =\left\| a\right\| ^3$ for all $a$ in ${A}.$
\end{enumerate}\smallskip

The following theorem is a celebrated result of I. Kaplansky (see
\cite[Theorem 6.2]{Kaps} or \cite[Theorem 1.2.4]{Sa}).

\begin{theorem}\label{t Kaplanski}
Let $A$ be a commutative C$^*$-algebra with a norm $\|.\|$ and let
$\|.\|_1$ be another norm on $A$ under which $A$ is a normed
algebra. Then $\|a\| \leq \|a\|_1$, for every $a$ in $A$. Further,
for any algebra norm, $\| .\|_{1}$, on $A_{sa}$, the inequality
$\|a\| \leq \|a\|_1$ holds for every $a$ in $A_{sa}$. $\hfill\Box$
\end{theorem}

Every C$^*$-algebra is a JB$^*$-triple with respect to the product
$\J abc = \frac12 \ ( a b^* c + cb^* a)$. It seems natural to ask
whether in the above Theorem \ref{t Kaplanski} the norm $\|.\|_1$
can be replaced with another norm $\|.\|_2$ under which $A$ is a
normed Jordan triple. The complex statement in the following
result was established by K. Bouhya and A. Fern{\'a}ndez L{\'o}pez in
\cite[Proposition 13]{BouFer}. A detailed proof is included here
for completeness reasons.

\begin{lemma}\label{l previous Kaplansky triple norm} Let $L\subset \mathbb{R}_0^{+}$ be a subset
of non-negative real numbers satisfying that $L\cup \{0\}$ is a
compact. Let $C_0(L)$ denote the Banach algebra of all real or
complex valued continuous functions on $L\cup \{0\}$ vanishing at
zero (equipped with the supremum norm $\|.\|_\infty$). Suppose
that $\|.\|_2$ is a $\|.\|_\infty$-continuous norm on $C_0(L)$
under which $C_0(L)$ is a normed Jordan triple. Then $\|.\|_2$ is
equivalent to an algebra norm on $C_0(L)$, and consequently
$\|.\|_{\infty}$ and $\|.\|_2$ are equivalent norms. More
concretely, writing $ M = \sup\{\|x\|_2 :  \|x\|_{\infty} \leq
1\}$ we have $ \|a\|_\infty \leq M N(C_0(L),\|.\|_2)\
 \|a\|_2,$ for all $a\in C_0(L)$.
\end{lemma}

\begin{proof} Since $\|.\|_2$ is $\|.\|_\infty$-continuous, there exists a
 positive $M$ such that $\|x\|_2\leq M \|x\|_{\infty}$,
for all $x\in C_0(L)$.\smallskip

When $L$ is compact $C_0(L)$ coincides with the C$^*$-algebra of
all complex-valued continuous functions on $L$ or with the
self-adjoint part of that C$^*$-algebra. Let $1$ denote the unit
element in $C(L)$. Take $a,b$ in $C(L)$. Applying that $\|.\|_2$
is a triple norm we have
$$\|a \ b \|_2= \|\J a1b \|_2 \leq N(C_0(L),\|.\|_2) \ \|a\|_2 \ \|1\|_2 \ \|b\|_2 $$
$$\leq N(C_0(L),\|.\|_2)\ M \ \|a\|_2 \
\|b\|_2.$$ This shows that $\|.\|_2$ is equivalent to $M
N(C_0(L),\|.\|_2)\ \|.\|_2$, and the latter is an algebra norm on
$C_0(L)$.\smallskip

Suppose that $L$ is not compact. Take $a$ and $b$ in $C_0(L)$. For
each natural $n$, let $p_n$, $a_n$ and $b_n$ be
the functions in $C_0(L)$ defined by $$a_n (t):=\left\{%
\begin{array}{ll}
    0, & \hbox{if $t\in [0, \frac{1}{2n}]\cap L$;} \\
    \hbox{affine}, & \hbox{if $t \in [\frac{1}{2n},\frac{1}{n}]\cap L$;} \\
    a(t), & \hbox{if $t\in [\frac{1}{n} , \infty)\cap L$.} \\
\end{array}%
\right. \ b_n (t):=\left\{%
\begin{array}{ll}
    0, & \hbox{if $t\in [0, \frac{1}{2n}]\cap L$;} \\
    \hbox{affine}, & \hbox{if $t \in [\frac{1}{2n},\frac{1}{n}]\cap L$;} \\
    b(t), & \hbox{if $t\in [\frac{1}{n} , \infty)\cap L$.} \\
\end{array}%
\right.
$$

$$\hbox{ and } p_n (t):=\left\{%
\begin{array}{ll}
    0, & \hbox{if $t\in [0, \frac{1}{4n}]\cap L$;} \\
    \hbox{affine}, & \hbox{if $t \in [\frac{1}{4n},\frac{1}{2n}]\cap L$;} \\
    1, & \hbox{if $t\in [\frac{1}{2n} , \infty)\cap L$.} \\
\end{array}%
\right.$$  Since $a_n \ b_n = \J {a_n}{p_n}{b_n}$ and
$\|p_n\|_{\infty}\leq 1,$  we deduce that
$$\|a_n \ b_n \|_2= \|\J {a_n}{p_n}{b_n} \|_2 \leq N(C_0(L),\|.\|_2)\ \|a_n\|_2 \ \|p_n\|_2 \ \|b_n\|_2 $$
$$\leq N(C_0(L),\|.\|_2)\ M
\ \|a_n\|_2 \ \|b_n\|_2.$$ Having in mind that $\|a_n-a\|_{\infty}
\to 0$, $\|b_n-b\|_{\infty} \to 0$, it follows, from the
$\|.\|_{\infty}$-continuity of the norm $\|.\|_2$, that
$$\|a \ b \|_2 \leq N(C_0(L),\|.\|_2)\ M \ \|a\|_2 \ \|b\|_2,$$ which shows
that $\|.\|_2$ is equivalent to $M N(C_0(L),\|.\|_2)\ \|.\|_2$,
and the latter is an algebra norm on $C_0(L)$. The final statement
is a direct consequence of Kaplansky's theorem (see Theorem \ref{t
Kaplanski}).
\end{proof}

\begin{remark}
{\rm Let $K$ be a compact Haussdorff space. Suppose that $\|.\|_2$
is a norm on $C(K)$ under which $C(K)$ is a normed Jordan triple
($\|.\|_\infty$-continuity of $\|.\|_2$ is not assumed). Let us
write $N=N(C(K),\|.\|_2)$. Following the argument given in the
proof of Lemma \ref{l previous Kaplansky triple norm}, we deduce
that
$$\|a \ b \|_2= \|\J a1b \|_2 \leq N\ \|1\|_2 \
\|a\|_2 \ \|b\|_2,$$ for all $a,b\in C(K)$, which shows that
$\|.\|_2$ is equivalent to $\|1\|_2\ N\ \|.\|_2$, and the latter
is an algebra norm on $C(K)$. It follows by Kaplansky's theorem,
that $\|a \|_{\infty} \leq \|1\|_2\ N\ \|a\|_2$, for all $a\in
C(K)$.}
\end{remark}

S.B. Cleveland applied Kaplansky's theorem to prove that every
continuous monomorphism from a C$^*$-algebra to a normed algebra
is bounded below (cf. \cite[Lemma 5.3]{CLev}), equivalently,
every C$^*$-algebra has MOANT. It follows as a consequence of
\cite[Theorem 1]{Ben} or \cite[Theorem
 10]{PeRiRo} or  \cite{HejNik}, that JB$^*$-algebras have MOANT.
In the setting of (complex) JB$^*$-triples, K. Bouhya and A. Fern{\'a}ndez L{\'o}pez
proved the following result:

\begin{proposition}\label{p Cleveland triples}\cite[Corollary
14]{BouFer} Let $T: E\to F$ be a continuous triple monomorphism
from a JB$^*$-triple to a normed complex Jordan triple. Then $T$ is
bounded below. That is, every JB$^*$-triple has MOTNT.
$\hfill\Box$
\end{proposition}

We shall complete this section by proving a generalization of
the above result to the setting of (real) J$^*$B-triples.\smallskip

We recall that a \emph{real JB$^*$-triple} is a norm-closed real
subtriple of a complex JB$^*$-triple (compare \cite{IsKaRo95}). A
\emph{J*B-triple} is a real Banach space $E$ equipped with a
structure of real Banach Jordan triple which satisfies $(JB^*2)$
and the following additional axioms:
\begin{enumerate} \item[{\small $(J^*B1)$}] $N({E})=1;$ \item[{\small
$(J^*B2)$}] $\sigma_{L(E)}^{\CC} (L(x,x)) \subset [0, + \infty)$
for all $x \in E$; \item[{\small $(J^*B3)$}] $\sigma_{L(E)}^{\CC}
(L(x,y)-L(y,x)) \subset i \RR$ for all $x,y \in E$.
\end{enumerate}

Every closed subtriple of a J$^*$B-triple is a J$^*$B-triple (cf.
\cite[Remark 1.5]{DaRu}). The class of J$^*$B-triples includes
real (and complex) C$^*$-algebras and real (and complex)
JB$^*$-triples. Moreover, in \cite[Proposition 1.4]{DaRu} it is
shown that complex JB$^*$-triples are precisely those complex
Jordan-Banach triples whose underlying real Banach space is a
J$^*$B-triple.\smallskip

T. Dang and B. Russo established a Gelfand theory for
J$^*$B-triples in \cite[Theorem 3.12]{DaRu}. This Gelfand theory
can be refined to show that given an element $a$ in a
J$^*$B-triple $E$, there exists a bounded set $L\subseteq
(0,\|a\|]$ with $L\cup \{0\}$ compact such that the smallest
(norm) closed subtriple of $E$ containing $a$, $E_a$, is
isometrically isomorphic to $$C_0(L,\RR) :=\{f\in C_0(L),
f(L)\subseteq \RR\},$$ (see \cite[Page 14]{BurPeRaRu}). The
argument given in the proof of Corollary 14 in \cite{BouFer} can
be adapted to prove the following result. The proof is included
here for completeness reasons.

\begin{proposition}\label{p Cleveland real triples} Let $T: E\to F$ be
a continuous triple monomorphism from a (real) J$^*$B-triple to a
normed Jordan triple. Then $T$ is bounded below. Equivalently,
every J$^*$B-triple has MOTNT.$\hfill\Box$.
\end{proposition}

\begin{proof} Take an arbitrary element $a$ in $E$. Let $E_a$
denote the smallest (norm) closed subtriple of $E$ containing $a$.
By \cite[Page 14]{BurPeRaRu}, there exists a subset $L\subseteq
(0,\|a\|]$ with $L\cup\{0\}$ compact satisfying that $E_a$ is
isometrically J$^*$B-triple isomorphic to $C_0(L,\RR),$ when the
latter is equipped with the supremum norm $\|.\|_{\infty}$. We
shall identify $E_a$ and $C_0(L,\RR)$. The mapping $T|_{E_a} : E_a
\cong C_0(L,\RR) \to F$ is a continuous triple monomorphism.
Therefore the mapping $x\mapsto \|x\|_2 :=\|T(x)\|$ defines a
$\|.\|_\infty$-continuous norm on $C_0(L,\RR)$ under which
$C_0(L,\RR)$ is a normed Jordan triple.\smallskip

Noticing that $N(E_a,\|.\|_2)\leq N(F)$ and $$M= \sup \{ \|x\|_2:
x\in E_a, \|x\|_\infty \leq 1\}\leq \|T\|,$$ Lemma \ref{l previous
Kaplansky triple norm} assures that $\|a\| \leq N(F) \ \|T\| \
\|T(a)\|,$ for every $a\in E$.
\end{proof}

We recall that a subspace $I$ of a normed Jordan triple $E$ is a
\emph{triple ideal} if $\{E,E,I\}+\{E,I,E\} \subseteq I.$ The
quotient of a normed Jordan triple by a closed triple ideal is a
normed Jordan triple. It is also known that the quotient of a
JB$^*$-triple (resp., a J$^*$B-triple) by a closed triple ideal is
a JB$^*$-triple (resp., a J$^*$B-triple) (compare
\cite{Ka}).\smallskip

Let $T: E\to F$ be a continuous triple homomorphism from a (real)
J$^*$B-triple to a normed Jordan triple. The kernel of $T$,
$ker(T)$, is a norm-closed triple ideal of $E$ and the linear
mapping $\widetilde{T}:E/ker(T) \rightarrow F$ given by $
\widetilde{T}(a +ker(T))=T(a)$ is a continuous triple monomorphism
from a (real) J$^*$B-triple to a normed Jordan triple and $
\widetilde{T}(E)=T(E).$ Proposition \ref{p Cleveland real triples}
assures that $\widetilde{T}$ is bounded bellow, and hence it has
closed range.\smallskip

A real JB$^*$-algebra is a closed $^*$-invariant real subalgebra
of a (complex) JB$^*$-algebra. Real C$^*$-algebras (i.e., closed
$^*$-invariant real subalgebras of C$^*$-algebras), equipped with
the Jordan product $a\circ b = \frac12 (a b +b a)$, are examples
of real JB$^*$-algebras.

\begin{corollary}\label{c cont closed range} Every continuous
triple homomorphism from a (real) J$^*$B-triple to a normed Jordan
triple has closed range. In particular, every continuous triple
homomorphism from a real or complex C$^*$-algebra to a normed
Jordan triple has closed range.$\hfill \Box$
\end{corollary}

\begin{corollary}\label{c real JB*-algebras} Let $A$ be a real
JB$^*$-algebra and let $B$ be a real Jordan Banach algebra (or a
real Jordan-Banach triple). Then every continuous triple
monomorphism from $A$ to $B$ is bounded below. That is, $A$ has
MOTNT and MOANT.$\hfill \Box$
\end{corollary}

\begin{corollary}\label{c real C*-algebras} Let $A$ be a real or complex
C$^*$-algebra and let $B$ be a real Banach algebra (or a real
Jordan-Banach triple). Then every continuous triple monomorphism
from $A$ to $B$ is bounded below. That is, $A$ has MOTNT and
MOANT.$\hfill \Box$
\end{corollary}

\section{Separating spaces for triple homomorphisms}

We have seen in the previous section that real and complex
C$^*$-algebras and real and complex JB$^*$-algebras have MOTNT and
MOANT. Equivalently, if $A$ denotes a real or complex
C$^*$-algebra (resp., a real or complex JB$^*$-algebra) every
continuous (triple) monomorphism $T$ from $A$ to a Banach algebra
(resp., a Jordan Banach algebra) is bounded below. C$^*$-algebras
and JB$^*$-algebras satisfy a stronger property: when $A$ is a
C$^*$-algebra (resp., a JB$^*$-algebra) every non-necessarily
continuous monomorphism from $A$ to a Banach algebra (resp., a
Jordan Banach algebra) is bounded below (compare \cite[Theorem
5.4]{CLev} and \cite[Theorem 1]{Ben} or \cite[Theorem
 10]{PeRiRo} or  \cite{HejNik}).\smallskip

The question clearly is whether every non-necessarily continuous
triple monomorphism from a complex JB$^*$-triple (resp., from a
real J$^*$B-triple) to a normed Jordan triple is bounded below. In
this section we provide a positive answer to this question.
Following a classical strategy, we shall study the
\emph{separating ideals} associated with a triple
homomorphism.\smallskip

Under additional geometric assumptions, triple homomorphisms are
automatically continuous. More concretely, every triple
homomorphism between two JB$^*$-triples is automatically
continuous (compare \cite[Lemma 1]{BarDanHor}). In this setting
the problem reduces to the question of minimality of triple norm
topology treated in section \S 2. However, when the codomain space
is not a JB$^*$-triple, the continuity of a triple homomorphism
does not follow automatically. We shall derive a new strategy to
solve Problem $(P)$ without any additional geometric hypothesis on
the codomain space.\smallskip

The following definitions and results are inspired by classical
ideas developed by C. Rickart \cite{Rick50}, B. Yood \cite{Yood},
W.G. Bade and P.C. Curtis \cite{BaCur} and S.B. Cleveland
\cite{CLev}.  Let $T: X\to Y$ be a linear mapping between two
normed spaces. Following \cite[Page 70]{Rick}, the
\emph{separating space}, $\sigma_Y (T)$, of $T$ in $Y$ is defined
as the set of all $z$ in $Y$ for which there exists a sequence
$(x_n) \subseteq X$ with $x_n \rightarrow 0$ and
$T(x_n)\rightarrow z$. The \emph{separating space}, $\sigma_X
(T)$, of $T$ in $X$ is defined by $\sigma_X (T):=T^{-1}(\sigma_Y
(T)).$ For each element $y$ in $Y$, $\Delta(y)$ is defined as the
infimum of the set $\{ \|x\|+\|y-T(x)\| : x\in X\}$. The mapping
$x\mapsto \Delta (x)$, called the \emph{separating} function of
$T$, satisfies the following properties:\begin{itemize}
\item[$a)$] $\Delta (y_1+y_2) \leq \Delta(y_1) + \Delta (y_2),$
\item[$b)$] $\Delta (\lambda y ) = |\lambda| \ \Delta (y)$,
\item[$c)$] $\Delta (y) \leq \|y\|$ and $\Delta (T(x))\leq \|x\|,$
\end{itemize}
for every $y,y_1$ and $y_2$ in $Y$, $x$ in $X$ and $\lambda$ scalar
(compare \cite[Page 71]{Rick} or \cite[Proposition 4.2]{CLev}).\smallskip

A straightforward application of the closed graph theorem shows
that a linear mapping $T$ between two Banach spaces $X$ and $Y$ is
continuous if and only if $\sigma_Y (T) =\{0\}$ (cf.
\cite[Proposition 4.5]{CLev}).\smallskip

It is not hard to see that $\sigma_Y (T) = \{ y\in Y: \Delta (y)
=0\}$, while $\sigma_X (T)= \{x\in X : \Delta (T(x)) = 0\}.$
Therefore $\sigma_X(T)$ and $\sigma_Y(T)$ are closed linear
subspaces of $X$ and $Y,$ respectively. The assignment $$x+
\sigma_X(T) \mapsto \widetilde{T} (x+ \sigma_X(T)) = T(x)
+\sigma_Y(T)$$ defines an injective linear operator from $X/
\sigma_X(T)$ to $Y/\sigma_Y(T)$. Moreover, $\widetilde{T}$ is
continuous whenever $X$ and $Y$ are Banach spaces.\smallskip

The separating subspaces of a triple homomorphism enjoy additional
algebraic structure.\smallskip

\begin{lemma}
\label{l separating spaces are ideals} Let $T: E \to F$ be a
non-necessarily continuous triple homomorphism between two normed Jordan triples.
Then $\sigma_E(T)$ is a norm-closed triple ideal of $E$ and
$\sigma_F(T)$ is a norm-closed triple ideal of the norm closure of
$T(E)$ in the completion of $F$.
\end{lemma}

\begin{proof}
Let us fix $z\in \sigma_F (T)$. In this case there exists a
sequence $(x_n) \subseteq E$ with $x_n \rightarrow 0$ and
$T(x_n)\rightarrow z$. Given $x,y$ in $E$, the sequences $(\J
{x_n}xy)$ and $(\J x{x_n}y)$ are norm-null, $$T(\J {x_n}xy) = \J
{T(x_n)}{T(x)}{T(y)} \to \J {z}{T(x)}{T(y)}$$ and $$ T(\J x{x_n}y)
= \J {T(x)}{T(x_n)}{T(y)} \to \J {T(x)}{z}{T(y)}.$$ This shows
that $\sigma_F(T)$ is a norm-closed triple ideal of
$\overline{T(E)}^{\|.\|}.$\smallskip

We have already proved that $$\J {\sigma_F (T)}{T(E)}{T(E)}
\subseteq \sigma_F (T)$$ and $$\J {T(E)}{\sigma_F (T)}{T(E)}
\subseteq \sigma_F (T).$$ This implies that $$T(\J {\sigma_E
(T)}{E}{E} ) \subseteq \J {\sigma_F (T)}{T(E)}{T(E)} \subseteq
\sigma_F (T)$$ and $$T(\J {E}{\sigma_E (T)}{E} ) \subseteq \J
{T(E)}{\sigma_F (T)}{T(E)} \subseteq \sigma_F (T),$$ which shows
that $\J {\sigma_E (T)}{E}{E}, \J {E}{\sigma_E (T)}{E} \subseteq
\sigma_E (T)$.
\end{proof}

The following result follows from Lemma \ref{l separating spaces
are ideals} and the basic properties of the separating spaces.

\begin{proposition}\label{p factor separating ideals} Let $T: E \to F$ be a
non-necessarily continuous triple homomorphism between two Jordan-Banach triples. Then the
mapping \linebreak $\widetilde{T} : E/\sigma_E(T) \rightarrow F/\sigma_F(T)$,
defined by $\widetilde{T}(a+E/\sigma_E(T))=T(a)+F/\sigma_F(T),$ is a continuous triple
monomorphism.$\hfill\Box$
\end{proposition}

Let $E$ be a normed Jordan triple. Two elements $a$ and $b$ in $E$ are said
to be \emph{orthogonal} (written $a\perp b$) if $L(a,b) = L(b,a)=0$. A direct application
of the Jordan identity yields that, for each $x$ in $E$,
\begin{equation}\label{eq orth inner ideal} a\perp \J bxb  \hbox{ whenever } a\perp b.\end{equation}

The following theorem is a ``\emph{main boundedness theorem}''
type result for Jordan-Banach triples (compare \cite[Theorem
2.1]{BaCur}, see also \cite[Theorem 3.1]{CLev}).

\begin{theorem}\label{t boundeness} Let $T:E\rightarrow F$ be a
non-necessarily continuous triple homomorphism
between Jordan-Banach triples and let $(x_n)$, $(y_n)$ be two
sequences of non-zero elements in $E$ such that $x_n \perp x_m, y_m$ for every
$n \neq m$, then
$$\sup \left\{ \frac{\| T(\{x_n,x_n,y_n \})\|}{\|x_n\|^2 \|y_n\|}, n \in \NN  \right\}<\infty .$$
\end{theorem}

\begin{proof} Suppose that
$\sup \left\{ \frac{\| T(\{x_n,x_n,y_n \})\|}{\|x_n\|^2 \|y_n\|},
n \in \NN  \right\}=\infty.$ Under this assumption, we may find a
subsequence $(a_{p,q})_{_{p,q\in \NN}}$ of $(x_n)$ formed by
mutually orthogonal elements such that
$$ \|T\{a_{p,q},a_{p,q},b_{p,q}\}\| > 4^p\ 8^q\ \|a_{p,q}\|^2 \|b_{p,q}\|,
\ \  p,q\in \NN,$$
 where $b_{p,q}=y_m$ whenever $a_{p,q}=x_m.$ Now, for each $p\in\NN$, we define $$z_p=\sum_{k=1}^\infty \frac{a_{p,k}}{2^k
 \|a_{p,k}\|}.$$ It is easy to see that, for each natural $q$, $b_{l,q}\perp z_p$
 whenever $l\neq p$. The equality $$\{z_p,z_p,b_{p,q}\}=\frac{1}{4^q \| a_{p,q}\|^2}\{a_{p,q},a_{p,q},b_{p,q}\},
 \  q \in \NN,$$ follows from the (joint) continuity of the triple
 product and the orthogonality hypothesis. Thus, $T(z_p)\neq 0,\ \forall p \in \NN$.\smallskip

For each $p$ in $\NN$ choose $n(p)$ in $\NN$ with $2^{n(p)}>\|
T(z_p)\|^2$ and define $y=\sum_{k=1}^\infty \frac{b_{k,n(k)}}{2^k
\|b_{k,n(k)}\|}.$ It follows that $$ \{ z_p,z_p,y\}=\frac{1}{2^{p}
4^{n(p)} \|b_{p,n(p)} \| \|a_{p,n(p)}\|^2}\
\{a_{p,n(p)},a_{p,n(p)},b_{p,n(p)}\}.$$ Therefore, $$N(F) \
\|T(y)\|\|T(z_p)\|^2>\|T(\{z_
p,z_p,y\})\|>2^{p}2^{n(p)}>2^p\|T(z_p)
 \|^2.$$ This implies that $N(F) \ \|T(y)\|>2^p$ for every positive integer $p,$
 which is impossible.
\end{proof}

Given an element $a$ in a normed Jordan triple $E$, we
denote $a^{[1]} = a$, $a^{[3]} = \J aaa$ and $a^{[2 n +1]} := \J
a{a^{[2n-1]}}a$ $(\forall n\in \mathbb{N})$. The Jordan identity implies
that $a^{[5]} = \J aa{a^{[3]}}$ , and by induction, $a^{[2n+1]} = L(a,a)^n (a)$
for all $n\in\NN$. The element $a$ is called \emph{nilpotent} if $a^{[2n+1]}=0$
for some $n$.\smallskip

A Jordan-Banach triple $E$ for which the vanishing of
$\J aaa$ implies that $a$ itself vanishes is said to be
\emph{anisotropic}. It is easy to check that $E$ is anisotropic
if and only if zero is the unique nilpotent element in $E$.\smallskip

Let $a$ and $b$ be two elements in an anisotropic normed Jordan
triple $E$. If $L(a,b)=0$, then, for each $x$ in $E$, the Jordan
identity implies that $$\J {L(b,a)x}{L(b,a)x}{L(b,a)x} = 0,$$ and
hence $L(b,a)=0$. Therefore $a\perp b$ if and only if
$L(a,b)=0$.\smallskip

In the setting of (complex) JB$^*$-triples, every element admits
3rd- and 5th- square roots. In fact, a continuous functional
calculus can be derived from the Gelfand representation for
abelian JB$^*$-triples (cf. \cite[\S 1]{Ka}). Let $a$ be an
element in a JB$^*$-triple $E$. Denoting by $E_a$ the
JB*-subtriple generated by the element $a$, it is known that $E_a$
is JB*-triple isomorphic (and hence isometric) to $C_0 (L)=
C_0(L,\CC)$ for some locally compact Hausdorff space $L\subseteq
(0,\|a\|],$ such that $L\cup \{0\}$ is compact. It is also known
that there exists a triple isomorphism  $\Psi$ from $E_a$ onto
$C_{0}(L)$ satisfying $\Psi (a) (t) = t$ $(t\in L)$ (compare
\cite[Lemma 1.14]{Ka} or \cite[Proposition 3.5]{Ka96}). Having in
mind this identification we can always find a (unique) element $z$
in $E_a$ such that $z^{[5]} = a$. The element $z$ will be denoted
by $a^{[\frac15]}$.\smallskip

When $E$ is a (real) J$^*$B-triple, we have already commented that
the norm closed subtriple generated by a single element $a$ is
triple isomorphic (and isometric) to $C_0(L,\RR) :=\{f\in C_0(L),
f(L)\subseteq \RR\},$ for some locally compact subset $L\subseteq
(0,\|a\|]$ with $L\cup \{0\}$ compact. Therefore there exists a
(unique) element $z=a^{[\frac15]}$ in $E_a$ such that $z^{[5]} =
a$.\smallskip

It should be noticed here that, in the setting of
J$^*$B-triples (resp., JB$^*$-triples) orthogonality is a ``local
concept'' (compare Lemma 1 in \cite{BurFerGarMarPe} whose proof
remains valid for J$^*$B-triples). Indeed, two elements $a$ and
$b$ in a J$^*$B-triple $E$ are orthogonal if and only if one of the
following equivalent statements holds: $$ (a) \ \J aab =0, \ \ (b)
\ E_a \perp E_b, \ \ (c) \ \{b,b,a\} =0,$$
$$(d)\  a\perp b \hbox{ in a subtriple of $E$ containing both
elements}.$$ It can be easily seen that $a\perp b$ if and only if $a^{[\frac15]} \perp b^{[\frac15]}$.

\begin{lemma}\label{l nonzero 5th power finite}
Let $T:E\rightarrow F$ be a
non-necessarily continuous triple homomorphism between two
Jordan-Banach triples and let $({x_n})$ be a sequence of mutually
orthogonal norm-one elements in $\sigma_E(T)$. Then, except for a
finite number of values of $n$, $T(x_n)^{[5]}=0.$ Further, if $E$
is a JB$^*$-triple or a (real) J*B-triple or $F$ is an anisotropic
Jordan Banach triple then $T(x_n)=0,$ except for finitely many $n
\in \NN.$
\end{lemma}

\begin{proof} We shall argue by contradiction, supposing that
$T(x_n)^{[5]}\neq 0$ for infinitely many $n$ in $\NN.$ By passing
to a subsequence, we may assume $T(x_n)^{[5]}\neq 0$ for every
$n\in \NN$. Since $(x_n)$ is a sequence in $\sigma_E(T),$ for each
$n\in \NN$, there is a sequence $(a_{n,k})_k\subseteq E$ such that
$\lim_k a_{n,k}=0$ and $\lim_k T(a_{n,k})=T(x_n).$ Thus, for each
$n$ in $\NN,$ $\lim_k \{x_n,a_{n,k},x_n\}=0.$ The continuity of
the triple product in $F$ implies that
$$ \lim_k T(\{x_n,x_n, \{x_n,a_{n,k},x_n \}\} )$$ $$=\lim_k \{
T(x_n),T(x_n),\{T(x_n),T(a_{n,k}),T(x_n) \} \}=T(x_n)^{[5]}\neq
0.$$

We observe that, for each $n\in\NN$, the set $$\{k\in \NN :
\{x_n,a_{n,k},x_n\}\neq 0\}$$ is infinite. Passing to a
subsequence of $(a_{n,k})$ we may assume that
$$\{x_n,a_{n,k},x_n\}\neq 0, \forall (n,k)\in
\NN\times \NN.$$

Therefore, $$ \lim_k \frac{\|T(\{x_n,x_n, \{x_n,a_{n,k},x_n \}\} )
\|}{\|\{x_n,a_{n,k},x_n\} \| }=\infty.$$

For each positive integer $n$, pick $m(n)$ such that
\begin{equation}\label{eq 2} \frac{\|T(\{x_n,x_n,
\{x_n,a_{n,m(n)},x_n \}\} ) \|}{\|\{x_n,a_{n,m(n)},x_n\} \| }>n\|x_n\|^2.\end{equation}

Writting $y_n=\{x_n,a_{n,m(n)},x_n \},$ it follows by $(\ref{eq orth inner ideal})$
that $y_n\perp x_m$ for $n\neq m$. The inequality $(\ref{eq 2})$ yields
$ \frac{\|T(x_n,x_n,y_n) \|}{\|x_n\|^2\|y_n\|}>n,
\forall n\in \NN,$ which contradicts the main boundeness theorem
(compare Theorem \ref{t  boundeness}).\smallskip

If $E$ is a JB$^*$-triple (resp., a J$^*$B-triple), by Lemma
\ref{l separating spaces are ideals}, $\sigma_E(T)$ is a norm
closed ideal of $E$ and hence a JB$^*$-triple (resp., a
J$^*$B-triple). Therefore, the sequence of mutually orthogonal
elements $(z_n)=(x_n^{_{_{[\frac{1}{5}]}}})$ lies in
$\sigma_E(T).$ Since $T(z_n)^{[5]} = T(z_n^{[5]}) = T(x_n)$, we
have $T(x_n)=0$ for finitely many $n$ in $\NN.$\smallskip

Finally, when $F$ is anisotropic the final statement follows
straightforwardly.\end{proof}

An element $e$ in a normed Jordan triple $E$ is called
\emph{tripotent} if $\J eee = e$. Every tripotent $e$ induces a
decomposition $E= E_{2} (e) \oplus E_{1} (e) \oplus E_0 (e)$ into
the corresponding \emph{Peirce spaces} where $E_j (e)$ is the
$\frac{j}{2}$ eigenspace of $L(e, e)$. Furthermore, the following
Peirce rules are satisfied,
$$ \J {E_{2} (e)}{E_{0}(e)}{E} = \J
{E_{0} (e)}{E_{2}(e)}{E} =0,$$
$$ \label{eq peirce rules2} \J
{E_{i}(e)}{E_{j} (e)}{E_{k} (e)}\subseteq E_{i-j+k} (e),$$ where
$E_{i-j+k} (e)=0$ whenever $i-j+k \notin \{ 0,1,2\}$ (compare
\cite[Proposition 21.9]{Up}). The projection $P_{j} (e): E \to E_j
(e)$ is called the \emph{Peirce-$j$} projection induced by
$e$.\smallskip

The Peirce-2 subspace, $E_2 (e)$, associated with a tripotent $e$
is a normed Jordan $^*$-algebra with respect to the product and
involution defined by \linebreak $x\circ_{e} y := \J xey$ and
$x^{\sharp_{e}} := \J exe$, respectively (compare, \cite[Lemma
21.11]{Up}).

\begin{lemma}\label{l trip in separating are sent to zero}
Let $T:E\rightarrow F$ be a
non-necessarily continuous triple homomorphism between two
Jordan-Banach triples. Then for each tripotent $e$ in $\sigma_E
(T)$ we have $T(e)= 0$.
\end{lemma}

\begin{proof} Suppose that there exists a tripotent $e$ in $\sigma_E(T)$
with $T(e)\neq 0.$ The linear mapping $T_{|E_2(e)}: E_2(e)
\rightarrow F_2(T(e))$ is a (unital) triple homomorphism between
(unital) Jordan Banach algebras. Then $T$ is a (unital) Jordan
homomorphism. Let $(x_n)$ be a sequence in $E$ such that $x_n
\rightarrow 0$ and $T(x_n) \rightarrow T(e).$ Then
$P_2(e)(x_n)\rightarrow 0$ and $T(P_2(e)(x_n))=P_2(T(e))(T(x_n))
\rightarrow T(e).$ Thus, $e$ is an idempotent in
$\sigma_{E_2(e)}(T_{|E_2(e)})$ with $T(e)\neq 0,$ which
contradicts Theorem 3.12 or Corollary 3.13 in
\cite{PuYood}.\end{proof}

\begin{lemma}\label{l monomorphs have zero separating space}
Let $T:E\rightarrow F$ be a
non-necessarily continuous triple monomorphism from a
JB$^*$-triple (resp., a J$^*$B-triple) to a Jordan-Banach triple. Then
$\sigma_E(T)=0.$
\end{lemma}

\begin{proof} Suppose that $\sigma_E(T)\neq 0.$ Then, by Lemma \ref{l separating spaces are ideals},
 $ \sigma_E(T)$ is a norm-closed triple ideal of $E,$ and hence a JB$^*$-triple
(resp., a J$^*$B-triple). Suppose that $a$ is a nonzero element in
$\sigma_E(T).$ We have already seen that, $E_a$  is isometrically
triple isomorphic to $C_0(L),$ for some subset $L\subseteq
(0,\|a\|]$ with $L\cup \{0\}$ compact.\smallskip

We claim that $L$ is finite. Indeed, if $L$ were infinite we could find, via
Urysohn's lemma, a sequence of mutually orthogonal norm-one elements
$(x_n)_{_{n \in \NN}}\subseteq E_a \subseteq \sigma_{E} (T)$. Since $T$ is injective
we have $T(x_n)\neq 0, \forall n\in \NN,$ which contradicts Lemma \ref{l
nonzero 5th power finite}. Therefore $L$ must be finite.\smallskip

Let $ t \in L.$ Since $L$ is finite, the function $e$ defined by
$e(t)=1$ and $e(L \setminus \{t\})=0$ lies in $C_0(L)$. The
element $e$ is a tripotent in $E_a\subseteq \sigma_E(T)$ with
$T(e)\neq 0,$ which, by Lemma \ref{l trip in separating are sent
to zero}, is impossible.\end{proof}

The following proposition is a direct consequence of Lemma \ref{l
monomorphs have zero separating space} and Proposition \ref{p
factor separating ideals}.

\begin{proposition}\label{p penultimate} Let $T:E\rightarrow F$ be a
non-necessarily continuous triple monomorphism from a (complex)
JB$^*$-triple (resp., a (real) J$^*$B-triple) to a Jordan-Banach triple.
Then the linear mapping $\widetilde{T}: E \rightarrow F/\sigma_F(T),$
$\widetilde{T}(a)=T(a)+F/\sigma_F(T),$
is a continuous triple monomorphism.$\hfill\Box$
\end{proposition}

\begin{theorem}\label{t kaplansky without conrinuity}  Let $T:E\rightarrow F$ be a
non-necessarily continuous triple monomorphism from a (complex)
JB$^*$-triple (resp., a (real) J$^*$B-triple) to a normed Jordan triple. Then
$T$ is bounded below.
\end{theorem}

\begin{proof} We may assume, without loss of generality, that  $F$ is
Jordan-Banach triple, otherwise we can replace $F$ with its canonical completion.\smallskip

Let $\pi $ denote the canonical projection of $F$ onto $F/\sigma_F(T)$.
Proposition \ref{p penultimate} assures that the linear mapping
$\widetilde{T}: E \rightarrow F/\sigma_F(T)$, $x\mapsto \pi (T(x))$,
is a continuous triple monomorphism. By Propositions \ref{p Cleveland real triples}
and \ref{p Cleveland triples}, there exists a positive constant $M$ satisfying
$$M \ \|x\| \leq \|\widetilde{T} (x)\| = \|\pi (T(x))\| \leq \|T(x)\|, \ (x\in E),$$
which shows that $T$ is bounded below.
\end{proof}

The following corollary is the desired generalisation of a result
due to B. Yood \cite{Yood54} and S.B. Cleveland \cite{CLev}.

\begin{corollary}\label{c Yood-Cleveland} Let $T:E\rightarrow F$ be a
non-necessarily continuous triple monomorphism from a (complex)
JB$^*$-triple (resp., a (real) J$^*$B-triple) to a normed Jordan
triple. Then the norm closure of $T(E)$ in the canonical
completion of $F$ decomposes as the direct sum of $T(E)$ and
$\sigma_{F} (T)$.
\end{corollary}

\begin{proof}
Let $b$ be an element in the norm closure of ${T(E)}$ in the
completion of $F$. By assumptions, there exists a sequence $(x_n)$
in $E$ such that $b = \lim T(x_n)$.\smallskip

Since, by Theorem \ref{t kaplansky without conrinuity}, $T$ is
bounded below, the sequence $(x_n)$ is a Cauchy sequence in $E$.
Therefore there exists $x_0$ in $E$ satisfying $\lim x_n - x_0 =
0$ and $\lim T(x_n -x_0) = b - T(x_0)$. This shows that $b =
T(x_0) + (b-T(x_0))$, where $b-T(x_0)\in \sigma_{F} (T)$. Finally,
$T(E) \cap \sigma_{F} (T) = T(\sigma_{E} (T)) =\{0\}$, by Lemma
\ref{l monomorphs have zero separating space}.
\end{proof}

\bigskip\bigskip

\end{document}